\documentclass[12pt]{amsart}
\usepackage{latexsym, amssymb, amsmath, amscd}
 \usepackage{eucal}

    \newtheorem{rema}{Remark}[section]

   \newtheorem{theo}[rema]{Theorem}
   \newtheorem{def-theo}[rema]{Definition-Theorem}

\newcommand{\Q}{{\mathbb{Q}}}
\newcommand{\Z}{{\mathbb{Z}}}

\newcommand{\uu}{{\vec{u}}}
\newcommand{\vv}{{\vec{v}}}

\newcommand{\vsp}{\vspace{1em}}

\newcommand{\spn}{{\mathrm{Span}}}

\newcommand{\prf}{{\em Proof.\,\,\,\,}}

\newcommand{\proj}{{\mathrm{Proj}}}

\newcommand{\spec}{{\mathrm{Spec}}}

\newcommand{\ds}{\displaystyle}

\newcommand{\ve}{\vec{e}}

\newcommand{\cA}{\mathcal{A}}

\renewcommand{\P}{\mbox{${\mathbb{P}}$}}

\newtheorem{thm}[rema]{Theorem}
\newtheorem{prop}[rema]{Proposition}
\newtheorem{lemma}[rema]{Lemma}
\newtheorem{cor}[rema]{Corollary}
\newtheorem{corollary}[rema]{Corollary}
\newtheorem{dfn}[rema]{Definition}
\newtheorem{example}[rema]{Example}

\newtheorem{remark}[rema]{Remark}

\title[A Family of Maximal Mathieu Subspaces of Matrix Algebras]
{A Family of Maximal Mathieu Subspaces of Matrix Algebras}
 
   \author{George F. Seelinger and Wenhua Zhao}      
    \date{\today}
 \address{G. F. Seelinger, Department of Mathematics, Illinois State University, Normal, IL 61761. {\it Email}: gfseeli@ilstu.edu}
\address{W. Zhao, Department of Mathematics, Illinois State University, Normal, IL 61761. {\it Email}: wzhao@ilstu.edu}
 
\date{\today}

\begin{document}

\begin{abstract}
Let $F$ be a field. In this note we give a construction for a family of maximal Mathieu subspaces (or Mathieu-Zhao subspaces) of the matrix algebras $M_n(F)$ $(n\ge 2)$. As an application we also give a classification of Mathieu subspaces of $M_2(F)$ under the condition that $F$ is algebraically closed.
\end{abstract}

   

 

\thanks{The second author has been partially supported by the Simons Foundation grant 278638.}

 \bibliographystyle{alpha}
    \maketitle


\renewcommand{\theequation}{\thesection.\arabic{equation}}
\renewcommand{\therema}{\thesection.\arabic{rema}}
\setcounter{equation}{0}
\setcounter{rema}{0}
\setcounter{section}{0}

\section{\bf Background and Motivation}

Let $F$ be a field. We start with the definition of a Mathieu subspace for an associative algebra $\mathcal A$ over $F$, which was introduced by the second author in \cite{GIC} (see also \cite{MS}).  

\begin{dfn}  A $F$-subspace $M$ of $\cA$ is a {\em Mathieu subspace} (MS) of $\cA$ if for all $a, b, c\in \cA$ such that $a^m\in M$ for all $m\geq 1$, there exists an $N$ (depending $a, b, c$) such that $b a^m c \in M$ for all $m \ge N$.
\end{dfn}

A MS (Mathieu subspace) is also called a Mathieu-Zhao space or Mathieu-Zhao subspace in the literature (e.g., see \cite{EKC, DEZ}). The introduction of MSs was directly motivated by the Mathieu Conjecture \cite{Ma} and the Image Conjecture \cite{IC}, each of which implies the Jacobian conjecture \cite{K, BCW, E}. For example, the Jacobian conjecture will follow if some explicitly given subspaces of multivariate polynomial algebras over $\mathbb C$ can be shown to be MSs of the polynomial algebras. For details, see \cite{IC, DEZ, EKC}. 

Note that ideals are MSs, but not conversely. Therefore  the concept of MSs can be viewed as a natural generalization of the concept of ideals. However, in contrast to ideals, MSs are currently far from being well-understood. This is even the case for the most of finite rings or finite dimensional algebras over a field. For example, the classification of all MSs of the matrix algebras $M_n(F)$ $(n\ge 3)$ is still wide open.

The study of MSs of matrix algebras $M_n(F)$ was initiated by the second author in \cite{MS} in which the following two results were proven.  

\begin{theo}{\cite[Thm.\ 4.2]{MS}}\label{ZIT}  Let $V\subseteq M_n(F)$ be a proper subspace of $M_n(F)$.  Then $V$ is a MS of $M_n(F)$ if and only if it does not contain any nonzero idempotent of $M_n(F)$.
\end{theo}

\begin{theo}\cite[Thm $5.1$]{MS}\label{MS5.1}
Let $H$ be the subspace of $M_n(F)$ consisting of all trace-zero matrices. Then the following statements holds:
\begin{enumerate}
  \item[$i)$] if char.\,$F=0$ or char.\,$F>n$, then $H$ is the only MS of $M_n(F)$ of codimension $1$;  
  \item[$ii)$] if char.\,$F\le n$, then $M_n(F)$ has no MS of  codimension $1$. 
\end{enumerate}
\end{theo}

Actually, the theorem above holds also for one-sided codimension one MSs of $M_n(F)$. 

It is easy to see from the two theorems above that every subspace of a MS of $M_n(F)$ is also a MS of $M_n(F)$. Therefore, to classify all MSs of $M_n(F)$ it suffices to classify all maximal MSs of $M_n(F)$. 

In particular, when char.\,$F=0$ or char.\,$F>n$, every subspace of the trace-zero codimension $1$ subspace $H$ is a MS of $M_n(F)$. 
Conversely, under the same condition on char.\,$F$ and $n\ge 3$, 
A. Konijnenberg \cite[Cor.\,3.5]{Kon} in his master thesis directed by A. van den Essen proved that every MS of $M_n(F)$ of codimension $2$ is a subspace contained in $H$. Furthermore, M. de Bondt \cite{deBondt} proved that it is also the case for all the MSs of $M_n(F)$ of codimension less than $n$.   

\begin{thm}{\cite[Thm.\,1.4]{deBondt}}\label{deBondtThm} 
Assume $\text{char.}\,F=0$ or $\text{char.}\,F\geq n$.  
If $V \subseteq M_n(F)$ is a MS of $M_n(F)$ of codimension less than $n$, then $V$ is a subspace of $H$.
\end{thm}
 
Therefore, we have a unique maximal MS of $M_n(F)$ of codimension $1$ (with the above restriction on the characteristic of $F$) and there are no maximal MSs of $M_n(F)$ with codimension between $2$ and $n-1$.  Note, in \cite[Prop.\ 1.2]{deBondt} a MS of $M_n(F)$ of codimension $n$, which is not contained in $H$, is given (hence is necessarily maximal) and de Bondt notes, based on \cite{Kon}, that ``{\it . . . the codimension $n$ case seems quite difficult}".

In Section \ref{S2} of this note we give a construction of a family of maximal MSs of $M_n(F)$, which for all $n\geq 2$ contains some maximal MSs of codimension $n$. In Section \ref{S3} we give a classification for all the (maximal) MSs 
of $M_2(F)$ under the condition that $F$ is algebraically closed.

\renewcommand{\theequation}{\thesection.\arabic{equation}}
\renewcommand{\therema}{\thesection.\arabic{rema}}
\setcounter{equation}{0}
\setcounter{rema}{0}

\section{\bf A Class of Maximal Mathieu Subspaces of $M_n(F)$}\label{S2}

%
Throughout this section $F$ stands for a field of arbitrary characteristic, whose algebraic closure is denoted by $\bar F$. 
We make use of the non-degenerate Killing form on $M_n(F)$ given by
\[ \langle b,c \rangle = \text{Tr}\,(bc) \mbox{ for all  } b,c\in M_n(F). \] 
Hence $H$, the subspace of trace zero matrices in $M_n(F)$, is equal to $I_n^\perp = \{ b\in M_n(F) : \text{Tr}\,(b) = 0\}$. 
More generally, for any set $S$, we
write $S^\perp = \{ b\in M_n(F) : \text{Tr}\,(bx) = 0 \mbox{ for all } x\in S\}$ for the orthogonal complement of $S$ 
with respect to the Killing form.

Now let $T\subseteq GL_n(F)$ be a (split, but not necessarily maximal) torus acting on $M_n(F)$ by conjugation and 
\[ Z(T,M_n): = \{ 
a \in M_n(F) : \tau a \tau^{-1} = a \mbox{ for all } \tau \in T \} \]
be the centralizer of $T$ in $M_n(F)$. Denote by 
 $X(T)$ the character group of $T$ and let $W(T,M_n)\subseteq X(T)$ 
be the set of weights of $T$ in $M_n(F)$. Since $X(T)$ is an abelian group, we will use $+$ as the binary operation of $X(T)$ and
denote the identity element by $0$.  For each $\omega\in W(T,M_n)$, let $M_\omega\subseteq M_n(F)$ be the weight space corresponding to $\omega$.   So $M_0 = Z(T,M_n)$.  Then we have the following direct sum composition: 
 
\[ M_n(F) = \bigoplus_{\omega\in W(T,M_n)} M_{\omega}. \]

For each $\omega\in W(T,M_n)$ we let $\pi_{\omega}:M_n(F)\rightarrow M_{\omega}$ be the corresponding projection.  We will use the observation that for any $a\in M_n(F)$ we have $\text{Tr}\,(a) = \text{Tr}\,(\pi_0(a))$.  
Finally, let $X_\Q = X(T)\otimes_{\Z}\Q$ and let $f:X_{\Q}\rightarrow \Q$ be a linear functional such that $f(\omega)\neq 0$ for
all $0\neq \omega \in W(T,M_n)$.  Define $W_f^+=\{ \omega \in W(T,M_n) : f(\omega) >0\}$ .

\begin{thm} \label{Class1} 
Let $\Lambda\in T$ be such that $\text{Tr}\,(\Lambda e) \neq 0$ for all nonzero idempotents $e\in M_0$. Let $V$ be a subspace of $M_n(F)$ such that 
both $V$ and $\pi_0(V)$ are contained in $\Lambda^{\perp}$.   
If $\pi_{\omega}(V)\pi_{-\omega'}(V) = 0$ for all $\omega, \omega'\in W_f^+$, then $V$ is a MS of $M_n(F)$.
\end{thm}

\prf  Assume that $V$ contains a nonzero idempotent $e\in V$.  Then $\pi_0(e) = \pi_0(e^2)$ and we can write
\[ \pi_0(e^2) = \pi_0(e)^2 + \sum_{\omega\in W_f^+} \left(\pi_\omega(e)\pi_{-\omega}(e) + \pi_{-\omega}(e)\pi_{\omega}(e)\right) . \]
Since by assumption $\pi_{\omega}(e)\pi_{-\omega'}(e) = 0$ for all $\omega,\omega'\in W_f^+$,  the above equation reduces to 

\begin{align}\label{Class1-peq1}
 \pi_0(e) = \pi_0(e^2) = \pi_0(e)^2 + \sum_{\omega\in W_f^+}  \pi_{-\omega}(e)\pi_{\omega}(e) . 
\end{align}

Set $c\!:=\sum_{\omega\in W_f^+} \pi_{-\omega}(e)\pi_{\omega}(e)$. Then using the above assumption again, 
we have  $c^2=0$. 
Therefore,  by eq.\,(\ref{Class1-peq1}), $\pi_0(e)$  must have only
$0$ and $1$ as eigenvalues in $\bar F$ with at least one nonzero eigenvalue. 
Hence there exists a nonzero idempotent $e'\in M_0$ such that $\pi_0(e) = e' + w$ for some nilpotent $w\in M_0$. 
Furthermore, since $\Lambda w$ is nilpotent, we have $\text{Tr}\,(\Lambda e) = \text{Tr}\,(\Lambda \pi_0(e)) = \text{Tr}\,(\Lambda e') \neq 0$. But this contradicts our assumption that $V\subseteq \Lambda^\perp$. Therefore  $V$ does not contain any nonzero idempotent and is a MS of $M_n(F)$ by Theorem \ref{ZIT}.  \qed

\begin{example} \label{t-BlockCase}   Let $e_1,\ldots, e_t\in M_n(F)$ be a set of orthogonal idempotents so that 
$I_n = e_1 + e_2 + \cdots + e_t$.  Let $T = F^\times e_1 + \cdots + F^\times e_t$ and let $n_i$ be the rank of $e_i$ for 
$1\leq i\leq t$.  If $\varepsilon_i \in X(T,M_n)$ is the character given by
$\varepsilon_i(\alpha_1e_1 + \cdots + \alpha_t e_t) = \alpha_i$ for all $\alpha_1e_2+\cdots + \alpha_te_t\in T$, then 
$\{ \varepsilon_i : 1\leq i\leq t\}$ forms a basis of $X(T,M_n)$ as a $\Z$-lattice and $W(T,M_n) = \{ \varepsilon_i - \varepsilon_j : 1\leq i,j\leq t \}$.  Let 
$f:X(T,M_n)\otimes_{\Z}\Q \rightarrow \Q$ be given by $f(\varepsilon_i) = i$ for all $1\leq i\leq t$. Let $\Lambda = \sum_{i=1}^t \sigma_i e_i\in T$  be such that $\sum_{i=1}^t \sigma_i k_i \ne 0$ for all integer $t$-tuples  
$\vec{0}\ne (k_1,k_2,\ldots,k_t) \in \prod_{i=1}^t [0,n_i]$.
Then by Theorem \ref{Class1} we see that the subspace 
\[ V \!:=  \Big(M_0\oplus \bigoplus_{\omega\in W_f^+} M_{\omega}\Big)\cap \Lambda^\perp  =
\left( M_0 \cap \Lambda^\perp \right) \oplus \bigoplus_{\omega\in W_f^+} M_{\omega}.
 \]
is a MS of $M_n(F)$. 
\end{example}

Before proceeding further we fix the following notation in the setting of Example \ref{t-BlockCase} above. For any $1\le i, j\le n$ 
we set $M_{i, j}=e_iM_n(F)e_j$. Then 
$M_n(F)=\bigoplus_{1\le i, j\le n}M_{i,j}$ and hence, 
each $a\in M_n(F)$ can be written as 
$a=\sum_{1\le i, j\le n}a_{ij}$ for some 
$a_{i j}\in M_{i,j}$. 
In this case we also write 
$a$ formally as a $t\times t$ matrix 
$(a_{i j})$. We call this matrix 
{\bf the block matrix form} of $a$.  
Then, in terms of the block matrix forms the MS $V$ in Example \ref{t-BlockCase} above is the subspace formed by all the matrices whose block matrix forms are lower triangular with the blocks $(a_{11}, a_{22}, ..., a_{tt})$ on the diagonal such that $a_{ii}\in M_{i,i}$ and $\sum_{i=1}^t \sigma_i \text{Tr}\,(a_{ii})=0$.

When $F$ is algebraically closed, a much more concrete way to look at the MS $V$ above is as follows. Up to conjugations we may assume that the idempotent $e_i$ $(1\le i\le t$) is the $n\times n$ matrix with  blocks $(0_{n_1\times n_1}, \cdots, 0_{n_{i-1}\times n_{i-1}}, I_{n_i\times n_i}, 0_{n_{i+1}\times n_{i+1}},\cdots , 0_{n_{t}\times n_{t}})$ on the diagonal and zero elsewhere. In this case the weight space $M_{i,j}=e_i M_n(F) e_j$ is canonically isomorphic to the space of $n_i\times n_j$ matrices over $F$. The MS $V$ in Example \ref{t-BlockCase}
is formed by all the matrices which are in the (usual) lower 
triangular block form with the matrices $(a'_{11}, a'_{22}, ..., a'_{tt})$ on the diagonal such that $a'_{ii}\in M_{n_i\times n_i}(F)$ and 
$\sum_{i=1}^t \sigma_i \text{Tr}\,(a'_{ii})=0$. 
%
 
In the rest of this section we will freely use 
the notations fixed in and after Example \ref{t-BlockCase}. 
We will show in Corollary \ref{cor2.6} that 
the MSs in Example \ref{t-BlockCase} when $t=2$ are actually maximal MSs of $M_n(F)$. However, when $t\ge 3$,  the MSs in Example \ref{t-BlockCase} are not maximal, which can be seen from the following example.

\begin{example} \label{CounterEx}  
Let $V$ be as in the Example \ref{t-BlockCase} with $t\ge 3$. 
Fix $u \in e_1 M_n(F)e_2$ and $w \in e_2 M_n(F)e_3$ such that $uw\ne 0$. 
Then let $U=F(u+w)+V$. We show below that $U$ does not contain any nonzero idempotent of $M_n(F)$. Hence $U$ by Theorem \ref{ZIT} is a MS of $M_n(F)$, which contains $V$ as a proper subspace.  

Assume that $U$ contains a nonzero idempotent $e$. Then $e=\alpha(u+w)+v$ for some $\alpha \in F$ and $v\in V$. Since $e_1Ue_3=0$, 
we have     
$0=e_1 e e_3=\alpha^2 uw$. 
Therefore, $\alpha=0$ and $e\in V$, which is a contradiction, since by Theorem \ref{ZIT} $V$ does not contain any nonzero idempotent. 
\end{example}

Next we consider the following family of MSs of $M_n(F)$, 
which coincides the case $t=2$ of Example \ref{t-BlockCase} 
when $e_1$ or $e_2=0$. 

\begin{example}\label{2blockex}  Let $e_1,e_2,e_3$ be as in Example \ref{t-BlockCase} with $t=3$, except we here do not assume 
that all $e_i\ne 0$. 
 Let $T$, $\varepsilon_i$ $(1\le i\le 3)$, and $X(T,M_n)$ be as in 
 Example \ref{t-BlockCase} with $t=3$.  
%
 Then $\{ \varepsilon_1,  \varepsilon_2,\varepsilon_3\}$ is a basis 
 of $X(T,M_n)$.  
Let  $T'\subseteq T$ be the subtorus defined by $T'=F^\times(e_1+e_2)+F^\times  e_3$ and let $\varepsilon'_1 = 
\varepsilon_1|_{T'}$.  We define $f:X(T',M_n)\otimes_{\Z}\Q\rightarrow \Q$ by setting $f(\varepsilon'_1)=1$ and
$f(\varepsilon_3) = 3$. Hence $W(T',M_n) = \{ \omega, 0, -\omega\}$ where $\omega = \varepsilon_3-\varepsilon'_1$.
Let $\Lambda = \sigma_1 (e_1+e_2) +  
\sigma_2 e_3\in T'$ be such that
$\sigma_1\neq \sigma_2$ and $k_1\sigma_1 + k_2\sigma_2 \ne 0$ for all integer ordered pairs $\vec{0}\ne (k_1, k_2)\in [0,n_1+n_2]\times [0,n_3]$. Then it follows from Theorem \ref{Class1} that 
\begin{align}\label{2blockex-eq1}
 V  &= \Big( Z(T',M_n)+e_1 M_n(F) e_3 + e_3 M_n(F)  e_2\Big) \cap \Lambda^\perp \nonumber \\ 
&=\left( Z(T',M_n)\cap \Lambda^\perp\right)+ e_1 M_n(F) e_3 + e_3 M_n(F)  e_2
\end{align}
is a MS of $M_n(F)$ of codimension $(n_1+n_2)n_3  + 1$.
\end{example}

The main result of this section is the following theorem. 

\begin{thm}\label{2blockmax} The MSs in Example \ref{2blockex} are maximal MSs of $M_n(F)$.
\end{thm}

\prf  Let $V$  be a MS of the type given in Eq.\,(\ref{2blockex-eq1}) and let $e_i, \varepsilon_i$ $(i=1,2,3)$ and $\Lambda = \sigma_1 (e_1+e_2) + \sigma_2 e_3$ be as in Example 
\ref{2blockex}.  Fix any $w \not\in V$.  We show below that there exists a nonzero idempotent $Q$ in $V+Fw$. Hence $V+Fw$ by Theorem \ref{ZIT} is not a MS of $M_n(F)$, and 
$V$ is a maximal MS of $M_n(F)$. 

Note that we may assume 
$w \in  Z(T',M_n)+e_3M_{n}(F)e_1 +e_2M_{n}(F)e_3$.  
  Furthermore, if $w\in  Z(T',M_n)$, then $w\not\in \Lambda^\perp$, and we may choose $Q=I_n$, for $I_n\in Z(T',M_n) \subseteq Fw+V$.  Therefore, we may assume further $w\not\in  Z(T',M_n)$.   
Let $w_0\in Z(T', M_n)$, $w_1 \in e_3M_n(F)e_1$ and    
$w_2\in e_2M_n(F)e_3$ such that $w=w_0+w_1+w_2$. 
Then  $w_1$ and $w_2$ are not both zero. 
We now divide the proof into two cases.  

{\bf Case 1:}  {\it Assume $w_0\in \Lambda^\perp$}. Then 
$w_0\in Z(T', M_n)\cap \Lambda^\perp\subseteq V$. 
Replacing $w$ by $w-w_0$ we can assume $w_0=0$ and $w=w_1+w_2$.  Without loss of generality, assume $w_1\neq 0$ is of rank $r$.  Then there exists $v\in e_1M_n(F)e_3\subseteq V$ such  that 
$w_1v$ and $vw_1$ 
are both idempotents of rank  $r$ with $w_1vw_1=w_1$ and $vw_1v=v$. 
The existence of $v$ follows from the following general fact: {\it for each $a\in  M_{m\times k}(F)$ of rank $r$, there exists  $b\in  M_{k\times m}(F)$ such that 
$ab$ and $ba$ are both idempotents of rank $r$, which satisfy the equations $aba=a$ and $bab=b$.}
Indeed, we can find $P\in GL_m(F)$ and $Q\in GL_k(F)$  such that $a=P\begin{pmatrix} I_r &0\\
0&0  \end{pmatrix}_{m\times k} Q$. Then we let $b=Q^{-1}\begin{pmatrix} I_r &0\\
0&0  \end{pmatrix}_{k\times m} P^{-1}$. 

Now we claim that, for any $\beta\neq -1$, we have that 
\[ Q \!:=  (1+\beta)^{-1}((e_3+w_2+\beta v)(e_3+w_1) + \beta(e_3-w_1v)) \]
is an idempotent with trace equal to $\text{rank}\,(e_3)=n_3$.  
To show $Q$ is an idempotent, we 
write $Q$ in the block matrix form as: 
\[ Q = \frac{1}{1+\beta} \left[ \begin{array}{ccc} 
\beta vw_1 & 0 & \beta v \\
w_2w_1 & 0 & w_2 \\ w_1 & 0 & (1+\beta)e_3-\beta w_1v \end{array}\right] \]
where the $(i,j)$-block represents an element in $e_iM_n(F)e_j$.  Then it is straightforward to check $Q^2=Q$ using the identities $vw_1v = v$ and $w_1vw_1 = w_1$.  

Next, note that  
$ \text{Tr}\,(\Lambda(Q-w)) = \text{Tr}\,(\Lambda Q) = 
\sigma_2 n_3+\frac{(\sigma_1-\sigma_2)\beta r}{\beta + 1}$. 
Note also that by the assumption on $\sigma_1$ and $\sigma_2$ we have $r\sigma_1+(n_3-r)\sigma_2\ne 0$, since $n_3-r \ge 0$ as $w_1\in e_3M_n(F)e_1$ so $r=\mathrm{rank}(w_1)\leq \mathrm{rank}(e_3) = n_3$.

Choose 
$\beta = -\frac{\sigma_2 n_3}{r\sigma_1+(n_3-r)\sigma_2}$. 
Since $\sigma_1\neq \sigma_2$, we get $\beta \neq -1$ and
$\text{Tr}\,(\Lambda(Q-w))=0$. Hence we have $Q-w\in V$ 
and $Q\in Fw+V$.

{\bf Case 2:}  {\it Assume $w_0\not\in \Lambda^\perp$}.  Then $Fw_0 + (\Lambda^\perp\cap  Z(T',M_n)) =  Z(T',M_n)$, whence there
exist $0\neq \alpha\in F$ and $x_0\in (\Lambda^\perp\cap  Z(T',M_n))\subseteq V$ such that $\alpha w_0 + x_0 = e_3$.  Set
\[ Q\!:= (\alpha w+ x_0 + \alpha^2 w_2w_1) = (e_3 + \alpha w_1 +\alpha w_2 + \alpha^2 w_2w_1)  \]
which in the block matrix form is
\[ Q = \left[ \begin{array}{ccc} 0 & 0 & 0 \\ \alpha^2 w_2w_1 & 0 & \alpha w_2 \\ \alpha w_1 & 0 & e_3\end{array}\right]. \]
Then it follows that $Q^2=Q$. 
Since $w_2w_1\in e_2 M_n(F)e_1 \subseteq \Lambda^\perp \cap  Z(T',M_n)\subseteq V$, we get $Q\in Fw + V$.  \qed
\vsp

Note that by letting $e_2=0$ in Example \ref{2blockex} and Theorem \ref{2blockmax} it is easy to check that we have also the following corollary.

\begin{corollary}\label{cor2.6}
 Let $e_1\in M_n(F)$ be a rank $r$ idempotent with $0<r<n$ and 
$\sigma_1\neq \sigma_2\in F$ such that
$k_1  \sigma_1  + k_2\sigma_2 \neq 0$ for all integer ordered pairs 
$\vec{0}\ne (k_1, k_2)\in [0, r]\times [0, n-r]$.  Set $e_2=I_n-e_1$. 
Then 
$$
V=\Big( \sigma_1 e_1+\sigma_2e_2 + e_1 M_n(F)e_3 \Big)^\perp,$$ 
which in the block matrix form is
 \[ V = \left( \begin{array}{cc} \sigma_1 e_1 & e_1 M_n(F) e_2 \\ 0&   \sigma_2 e_2 \end{array}\right)^\perp \]
is a maximal MS of $M_n(F)$ of codimension $(rn-r^2+1)$.   
\end{corollary}
%

When $r=1$ the MSs in the corollary above have codimension $n$. 
The example given in \cite[Prop. 1.2]{deBondt} is 
one of these MSs of this type.

\renewcommand{\theequation}{\thesection.\arabic{equation}}
\renewcommand{\therema}{\thesection.\arabic{rema}}
\setcounter{equation}{0}
\setcounter{rema}{0}

\section{\bf Classification of Mathieu Subspaces of $M_2(F)$}\label{S3}

Throughout this section $F$ stands for 
an algebraically closed field of arbitrary characteristic. 
In this section we give a classification for all MSs 
of $M_2(F)$. As pointed out in Section $1$, by Theorems \ref{ZIT} and \ref{MS5.1} we only need to classify all maximal MSs of $M_2(F)$.
We start with the following lemma.

\begin{lemma}\label{lma3.1}
Let $V$ be a MS of $M_2(F)$ and $a\in V$ with $\text{Tr}\,(a)\ne 0$. 
Then the following two statements hold:
\begin{enumerate}
  \item[$i)$] $a$ is invertible;
  \item[$ii)$] if $\dim_F V=2$, then $a$ has distinct (nonzero) eigenvalues. 
 \end{enumerate}
\end{lemma}
\prf $i)$ Assume that $a$ is singular. Then $a$ is of rank $1$ and 
$a^2=\text{Tr}\,(a)a$ by the Cayley-Hamilton Theorem. Hence $\text{Tr}\,(a)^{-1}a$ is a nonzero idempotent in $V$, and,  by Theorem \ref{ZIT}, $V$ is not a MS of $M_2(F)$, contradicting our assumption on $V$.   

$ii)$ Let $\lambda_1$ and $\lambda_2$ be the eigenvalues of $a$. 
Then   $\lambda_1\lambda_2\ne 0$ by $i)$.  
Let $c\in I_2^\perp$ so that $\{a,c\}$ is a basis of $V$.  Note that for all $x\in F$ we have 
$\text{Tr}\,(a+xc)=\text{Tr}\,(a)\neq 0$.   
Then by $i)$ $\det(a+xc)\neq 0$ for all $x\in F$.  Since $F$ is algebraically closed, this is true if and only if 
$a^{-1}c$ has no nonzero eigenvalue in $F$.
Therefore, $a^{-1}c$ must be nilpotent, which gives us 
$a^{-1}ca^{-1}c = 0 \Rightarrow ca^{-1}c = 0$. Hence, 
$\det c=0$, so $c$ is nilpotent of degree $2$ by the Cayley-
Hamilton Theorem since we have $\text{Tr}\,(c) = 0$ as well. 
Viewing $c$ as a nonzero linear endomorphism of $F^2$ we have 
$\ker(c) = {\mathrm{Im}}(c)=\ker(ca^{-1})$ (as one-dimensional subspaces  of $F^2$).  Fix any $0\neq \vv \in \ker(c)$. Then 
$\vv$ is simultaneously an eigenvector of both $a$ and $c$ since $a\vv\in \ker(ca^{-1})=\ker(c)$. 
Without loss of generality, assume $a\vv = \lambda_1\vv$.

Assume $\lambda_1=\lambda_2$. Note that 
$a-\lambda_1I_2 \ne 0$ (otherwise $I_2$ would be in $V$)  
and is nilpotent. Therefore we have 
$$
\ker(a-\lambda_1I_2) = {\mathrm{Im}}(a-\lambda_1I_2) = F\vv = {\mathrm{Im}}(c) =
\ker(c).$$
Hence there exists $\beta\in F$ such that $a-\lambda_1 I_2=\beta c$. But this implies $I_2\in V$, which by Theorem \ref{ZIT} is a contradiction. Therefore $\lambda_1\ne\lambda_2$ and statment $ii)$ holds.
\qed \\

In his master thesis (\cite[Thm.\ 3.10]{Kon}) directed by A. van den Essen, A. Konijnenberg characterized all the codimension $2$ subspaces of $M_2(F)$ when $F$ is algebraically closed and 
$\text{char.}\,F \neq 2$.  
Now we can give a different proof for this result without the condition $\text{char.}\,F \neq 2$.

\begin{cor} \label{m2kprop} Let 
$V \subseteq M_2(F)$ be a MS with $\dim_F\, V =2$.  
Then either $V \subseteq I_2^\perp$ or there exist nonzero 
idempotents $e_1,e_2 \in M_2(F)$ such that $e_1+e_2 = I_2$ 
and $V=(\lambda_1e_1 + \lambda_2e_2)+ e_1M_2(F)e_2$ for some distinct nonzero $\lambda_1, \lambda_2\in F$ with  
$\lambda_1+\lambda_2\ne 0$. 

In other words, either $V \subseteq I_2^\perp$ or 
$V$ is conjugate to the subspace 
$\left\{ \begin{pmatrix} \lambda_1s & t \\
0 &\lambda_2s\end{pmatrix}\,\big|\, s, t\in F\right\}$ for some distinct nonzero $\lambda_1, \lambda_2\in F$ with  
$\lambda_1+\lambda_2\ne 0$.
\end{cor}

\prf Assume $V\not\subseteq I_2^\perp$.  Then there exists an $a\in V$ with $\text{Tr}\,(a)\neq 0$.  Let $\lambda_1$ and $\lambda_2$ be eigenvalues of $a$ in $F$.     
Then $\lambda_1+\lambda_2\ne 0$, and  
by Lemma \ref{lma3.1}, we also have 
$\lambda_1\lambda_2\ne 0$ and $\lambda_1\neq \lambda_2$.

Let $c$ and $\vv$ be as in the proof above for Lemma \ref{lma3.1}, $ii)$.
Then by the arguments there we have $ca^{-1}c=0$; $c\vv=0$ 
and $a\vv=\lambda_1v$. Let 
$\uu\in F^2$ be the eigenvector of $a$ such that $a\uu = \lambda_2\uu$.  Let $e_1, e_2\in M_2(F)$ be the idempotents
corresponding to the decomposition $F^2 = F\vv \oplus F\uu$.  
Then it is easy to see that $V$ is of the desired form.  
\qed  

\begin{remark}\label{rmk3.3} 
$i)$ For the subspace $V$ in Corollary \ref{m2kprop} such that 
$V\not \subseteq I_2^\perp$ we let $\sigma_1 = -\lambda_2$ and 
$\sigma_2 = \lambda_1$. Then $V$ is a subspace of the type given  
in Corollary \ref{cor2.6} by which  $V$ is indeed a MS of $M_2(F)$.  
 
$ii)$ If char.\,$F\ne 2$, then $I_2^\perp$ does not 
contain any nonzero idempotent. Hence each subspace 
$V\subseteq I_2^\perp$ by Theorem \ref{ZIT} is a MS of $M_2(F)$. 
Therefore, Corollary \ref{m2kprop} in this case coincides 
with \cite[Thm.\ 3.10]{Kon} 
which characterizes all MSs of codimension $2$ of $M_2(F)$.
\end{remark}

\begin{cor} \label{cor3.4}
$i)$ If char.\,$F\ne 2$, then the dimensional  
$1$ maximal MSs of $M_2(F)$ are exactly the subspaces 
of the form $F(I_2+c)$, where $c \in M_2(F)$ is nilpotent. 

$ii)$ If char.\,$F=2$, then there is no dimension $1$ maximal MS in 
$M_2(F)$. 
\end{cor}

\prf Let $V$ be a dimensional $1$ maximal MS of $M_2(F)$ with 
$V=Fa$ for some $a\in M_2(F)$.  
Then by Theorem \ref{ZIT} $a^2\ne \beta a$ 
for any nonzero $\beta\in F$, otherwise $\beta^{-1}a$ would be 
a nonzero idempotent in $V$.  
If $\text{Tr}\,(a)=0$, then $V\subseteq I_2^\perp$ 
which would obviously not be maximal. 
Therefore $\text{Tr}\,(a)\ne 0$ and, 
 by Lemma \ref{lma3.1}, $a$ is invertible. 

Let $\lambda_1, \lambda_2$ be the nonzero eigenvalues 
of $a$. Then $\lambda_1+\lambda_2\ne 0$. If $\lambda_1\ne \lambda_2$, then $V$ is a proper subspace of 
$(\lambda_1e_1 + \lambda_2e_2)+ e_1M_2(F)e_2$ in 
Corollary \ref{m2kprop}, which as pointed out in Remark \ref{rmk3.3}, $i)$ is a MS of $M_2(F)$, whence $V$ would not be maximal.   
Therefore we have $\lambda_1=\lambda_2$, 
which we may assume are both equal to $1$.  
Let $c=a-I_2$. Then $c$ is nonzero and nilpotent, 
and $V=F(I_2+c)$. 

Conversely, for any nonzero nilpotent $c\in M_2(F)$, 
$I_2+c$ is invertible. 
The subspace $F(I_2+c)$ 
does not contain any nonzero idempotents, 
as the only invertible idempotent of 
$M_2(F)$ is $I_2$. Then by Theorem \ref{ZIT} 
$F(I_2+c)$ is a MS of $M_2(F)$.

Now assume char.\,$F\ne 2$. To show statement $i)$, 
we only need to show that, for any nonzero nilpotent 
$c\in M_2(F)$, the subspace $F(I_2+c)$ 
is maximal among MSs of $M_2(F)$. 
Assume otherwise. Then $V$ is strictly contained in a proper 
MS $W$ of $M_2(F)$. If $\dim_F W=3$, then $W=I_2^\perp$ 
by Theorem \ref{MS5.1}, which is a contradiction of 
$\text{Tr}\,(I_2+c)=2\ne 0$. Therefore $\dim_F W=2$. 
Then by Lemma \ref{lma3.1} the eigenvalues of 
$I_2+c$ must be distinct, which is a contradiction. 
Therefore $F(I_2+c)$ is a maximal MS of $M_2(F)$, and 
statement $i)$ follows.

Now assume char.\,$F=2$. To show statement $ii)$ 
it suffices to show that, for all nonzero nilpotent 
$c\in M_2(F)$, the subspace $F(I_2+c)$ is not maximal 
among MSs of $M_2(F)$. Note that 
$\dim_F I_2^\perp=3$ and the only nonzero idempotent    
contained in $I_2^\perp$ is $I_2$. Hence there exists 
a dimensional $2$ subspace $U$ of $I_2^\perp$ such that  
$I_2+c \in U$ and $I_2\not\in U$. Then by Theorem \ref{ZIT},
$U$ is a MS of $M_2(F)$ which strictly contains $V$. 
Therefore $V$ is not a maximal MS of $M_2(F)$, 
and statement $ii)$ follows.  
\qed \\

Now we give a classification for maximal MSs of 
$M_2(F)$ in the next two theorems. 

\begin{thm}\label{2max} 
Assume char.\,$F\ne 2$. 
Then the maximal MSs of $M_2(F)$ are exactly the followings:
\begin{enumerate}
\item[$i)$] 
$I_2^\perp$; 
\item[$ii)$] the subspaces of the form $(\lambda_1 e_1 + \lambda_2 e_2)+ e_1M_2(F)e_2$  for some nonzero idempotents $e_1, e_2$ such that $e_1+e_2=I_2$ and some distinct nonzero $\lambda_1,\lambda_2\in F$ with $\lambda_1 + \lambda_2\neq 0$; 
\item[$iii)$] the subspaces of the form $F(I_2+c)$ for 
some nonzero nilpotent $c \in M_2(F)$. 
\end{enumerate}
\end{thm}
 
\prf By Theorem \ref{ZIT} the only MS of dimension $3$ of $M_2(F)$ 
is $I_2^\perp$. Consequently, all dimensional $2$ MSs 
$V$ of $M_2(F)$ with $V\not\subseteq I_2^\perp$ are 
maximal. Then, by Corollary \ref{m2kprop} and Remark \ref{rmk3.3}, $i)$ we see that  maximal MSs of $M_2(F)$ of dimension $2$ are exactly 
those given in $ii)$. Finally, by Corollary \ref{cor3.4} the maximal MSs of $M_2(F)$ of dimension $1$ are exactly those given in $iii)$.   
\qed

\begin{thm}\label{2.2max}
Assume char.\,$F=2$. Then the maximal MSs of $M_2(F)$ are exactly the followings: 
\begin{enumerate}
\item[$i)$] 
the dimensional $2$ subspaces $V$ of $I_2^\perp$ such that $I_2\not\in V$; 
\item[$ii)$]  
the subspaces of the form $(\lambda_1 e_1 + \lambda_2 e_2)+ e_1M_2(F)e_2$  for some nonzero idempotents $e_1, e_2$ such that $e_1+e_2=I_2$ and some distinct nonzero $\lambda_1,\lambda_2\in F$ such that $\lambda_1 + \lambda_2\neq 0$. 
\end{enumerate} 
\end{thm}

\prf First, since the only nonzero idempotent contained in 
$I_2^\perp$ is $I_2$, the subspaces given in $i)$ are MSs of $M_2(F)$.
By Remark \ref{rmk3.3}, $i)$ the subspaces given in $ii)$ are also MSs of $M_2(F)$. Since $M_2(F)$ by Theorem \ref{MS5.1} has no 
MS of dimension $3$, all MSs of dimension $2$ are maximal. 
Therefore the theorem follows immediately from 
Corollaries \ref{m2kprop}  
and \ref{cor3.4}.  
  \qed \\

We end this section with the following example to show that 
Theorems \ref{2max} and \ref{2.2max} do not hold if the base field 
$F$ is not algebraically closed. Furthermore, the example 
  shows also that in general MSs are  not preserved 
under base field extensions.

\begin{example}  
Let $K$ be a field such that there exists $s \in K$ with  
$s \ne 0, \pm 1$ and $s\ne t^2$ for all $t\in K$. 
Let $a=\begin{pmatrix}1 & 0 \\0& s \end{pmatrix}$, 
$b=\begin{pmatrix}0&1\\1&0 \end{pmatrix}$ and 
$V=Ka+Kb$. Since $\det(a+xb)=s-x^2\ne 0$ for all $x\in K$, 
all nonzero elements of $V$ are invertible. Hence 
$V$ by Theorem \ref{ZIT} is a MS of $M_2(K)$  
(for $I_2\not\in V$), and by Theorem \ref{MS5.1} is maximal (for $V\not\subseteq I_2^\perp$). 
However, $V$ cannot be any of MSs 
of dimension $2$  given in Theorem \ref{2max} or Theorem \ref{2.2max}, 
since each of those MSs that is not contained in $I_2^\perp$   contains some nonzero nilpotent elements but $V$ does not.   

Furthermore, let $L$ be a field that contains $K$ and $\sqrt{s}$.  
Set $U\!:=L \otimes_K V$ and $c\!:=\frac1{1+s} (a+\sqrt{s}\, b)$. 
Then $c$ is a nonzero idempotent in $U$, whence $U$ by 
Theorem \ref{ZIT} is not a MS of $M_2(L)=L\otimes_K M_2(K)$.  
\end{example}


\begin{thebibliography}{FLM2}
\bibitem[BCW]{BCW} H. Bass, E. Connell and D. Wright, {\it The Jacobian Conjecture, Reduction of Degree and Formal Expansion of the Inverse}. Bull.  Amer. Math.  Soc.  \textbf{7}, (1982), 287--330. 

\bibitem[Bo]{deBondt} M. de Bondt, {\it Mathieu Subspaces of Codimension Less Than $n$ of
${\mathrm Mat}_n(K)$}.  {\em Linear Multilinear Algebra}, {\bf 64}(10) (2016), 2049--2067.

\bibitem[DEZ]{DEZ} H. Derksen, A. van den Essen and W. Zhao, {\it The Gaussian Moments Conjecture and the Jacobian Conjecture}.  Israel J. Math. {\bf 219} (2017), no.\,2, 917--928. See also arXiv:1506.05192 [math.AC]. 

\bibitem[DK]{DK}  J. J. Duistermaat and W. van der Kallen, {\it Constant Terms in Powers of a Laurent Polynomial}. Indag. Math. (N.S.) {\bf 9}\, (1998), no. 2, 221--231.

\bibitem[E]{E} A. van den Essen, {\it Polynomial Automorphisms and the Jacobian Conjecture}. Prog. Math., Vol.190, Birkh\"auser Verlag, Basel, 2000. 

\bibitem[EKC]{EKC} A. van den Essen, S. Kuroda and A. J. Crachiola, {\it Polynomial Automorphisms
and the Jacobian Conjecture: New Results from the Beginning of the 21st
Century}. Frontiers in Mathematics. Birkh\"auser $2021$. 
 
\bibitem[Ke]{K} O. H. Keller, 
{\it Ganze Gremona-Transformationen}. Monats. Math. Physik {\bf 47} (1939), no.\,1, 299-306. 

\bibitem[Kon]{Kon} A. Konijnenberg, {\it Mathieu Subspaces of Finite Products of Matrix Rings} [Master thesis]. 
Radboud University of Nijmegen; 2012. Directed by A.R.P. van den Essen.

\bibitem[Mat]{Ma} O. Mathieu, {\it Some Conjectures about Invariant Theory and Their Applications.} Alg\`ebre non commutative, groupes quantiques et invariants (Reims, 1995), 263--279, S\'emin. Congr., 2, Soc. Math. France, Paris, 1997. 

\bibitem[Z1]{IC} W. Zhao, {\it Images of Commuting  Differential Operators of Order One with Constant Leading Coefficients}.  J. Alg. {\bf 324} (2010),  no. 2, 231--247. [MR2651354]. See also arXiv:0902.0210 [math.CV]. 

\bibitem[Z2]{GIC} W. Zhao, {\em Generalizations of the Image Conjecture and the Mathieu Conjecture}.  J. Pure Appl. Alg. {\bf 214} (2010), 1200-1216. See also arXiv:0902.0212 [math.CV].  
 
\bibitem[Z3]{MS} W. Zhao, {\em Mathieu Subspaces of Associative Algebras}. J. Alg. {\bf 350} (2012), no.2, 245-272. 
 

\end{thebibliography}
\end{document}